\documentclass[11pt]{article}
\usepackage{authblk}

\usepackage{geometry}
\usepackage{amsmath,amssymb,mathrsfs}
\usepackage{stmaryrd}
\usepackage{color}
\usepackage[dvipsnames]{xcolor}
\usepackage{graphicx}
\usepackage{subfigure,diagbox}
\usepackage{soul}
\usepackage{algorithm,algpseudocode}

\usepackage{lineno,hyperref}
\modulolinenumbers[5]

\DeclareMathOperator*{\Loss}{Loss}

\def\bx{\mathbf{x}}

\def\bn{\mathbf{n}}
\def\bb{\mathbf{b}}
\def\bv{\mathbf{v}}
\def\bft{\mathbf{\theta}}
\def\grad{\nabla}
\def\pd{\partial}
\newcommand{\beq}{\begin{equation}}
\newcommand{\eeq}{\end{equation}}
\newcommand{\beqs}{\begin{eqnarray}}
\newcommand{\eeqs}{\end{eqnarray}}
\newcommand{\beqsn}{\begin{eqnarray*}}
\newcommand{\eeqsn}{\end{eqnarray*}}
\newcommand{\bary}{\begin{array}}
\newcommand{\eary}{\end{array}}

\newcommand*{\abs}[1]{\left|#1\right|}
\newcommand*{\dbblk}[1]{\llbracket #1 \rrbracket}

\title{A cusp-capturing PINN for elliptic interface problems}

\author[1]{Yu-Hau Tseng}
\author[2, 4]{Te-Sheng Lin}
\author[3, 4]{Wei-Fan Hu}
\author[2]{Ming-Chih Lai}
\affil[1]{Department of Applied Mathematics, National University of Kaohsiung, Kaohsiung 81148, Taiwan}
\affil[2]{Department of Applied Mathematics, National Yang Ming Chiao Tung University, Hsinchu 30010, Taiwan}
\affil[3]{Department of Mathematics, National Central University, Taoyuan 32001, Taiwan}
\affil[4]{National Center for Theoretical Sciences, National Taiwan University, Taipei 10617, Taiwan}

\begin{document}

\maketitle

\begin{abstract}
In this paper, we propose a cusp-capturing physics-informed neural network (PINN) to solve discontinuous-coefficient elliptic interface problems whose solution is continuous but has discontinuous first derivatives on the interface. To find such a solution using neural network representation, we introduce a cusp-enforced level set function as an additional feature input to the network to retain the inherent solution properties; that is,  capturing the solution cusps (where the derivatives are discontinuous) sharply. In addition, the proposed neural network has the advantage of being mesh-free, so it can easily handle problems in irregular domains. We train the network using the physics-informed framework in which the loss function comprises the residual of the differential equation together with certain interface and boundary conditions. We conduct a series of numerical experiments to demonstrate the effectiveness of the cusp-capturing technique and the accuracy of the present network model. Numerical results show that even using a one-hidden-layer (shallow) network with a moderate number of neurons and sufficient training data points, the present network model can achieve prediction accuracy comparable with traditional methods. Besides, if the solution is discontinuous across the interface, we can simply incorporate an additional supervised learning task for solution jump approximation into the present network without much difficulty.
\end{abstract}


\section{Introduction}
\label{sec:intro}

The study of fluid-structure interaction~(FSI) problems has been an important research topic in fluid dynamics for centuries, with applications ranging from, for example, fundamental physics, engineering, geophysics, and biomedicine. Typical small-scale examples include collisions between droplets in interfacial flows~\cite{PTCHWL16,TB05}, the dynamics of red blood cells flowing in pulsating arteries~\cite{Ku97,SBS02}, and the electrophoretic motion of colloidal particles in electrically charged fluids~\cite{HSR03,OW78}. The key components in these examples are fluid flow, deformable interfaces, and the complex mechanisms behind them. Moreover, physical parameters (such as viscosity or density) for each subregion of the domain may be different, resulting in lower regularity of the solution across the interfaces, thus requiring additional treatments for accurate simulations.

For instance, when the no-slip boundary condition is applied to a fluid-structure interface, the velocity field in the FSI problem is continuous in the entire domain, but its derivative is discontinuous across the interface. Among many classical numerical methods for solving such problems, Peskin proposed the immersed boundary~(IB) formulation~\cite{Peskin77,Peskin02}, which transforms the core of solving the velocity field into an elliptic problem with singular forces. The IB method adopts a regularized version of the Dirac delta function to discretize the singular forces directly, resulting in only first-order solution accuracy~\cite{LM12}. Another way to write the velocity equations is to impose jump conditions directly on the interface. So the problem becomes an elliptic interface problem in which the solution is continuous, but its normal derivative has jump discontinuity across the interface, which is exactly the formulation we aim to solve in this work.

Since the introduction of the IB formulation, several jump-capturing and high-order methods have been proposed for elliptic interface problems with discontinuous coefficients. For instance, LeVeque and Li introduced the immersed interface method (IIM)~\cite{LL94}, incorporating the jump conditions via local coordinates into the finite difference scheme to achieve the overall second-order accuracy in maximum norm. A simple implementation version of IIM that directly uses the jump conditions without introducing local coordinates was developed in \cite{HLY15,LT08} to achieve second-order accuracy in maximum norm as well. Liu et al.~\cite{LFK00} introduced a boundary condition capturing method (also known as the ghost fluid method (GFM)) that is able to solve the elliptic interface problems in a dimension-by-dimension manner, and can capture the solution and its normal derivative jumps sharply. However, the original GFM smoothes its tangential derivative, so the method is only first-order accurate in the maximum norm. Egan and Gibou~\cite{EG20} extended the original GFM by recovering the convergence of the gradients to achieve second-order accuracy without modifying the resultant linear system. There are many other Cartesian grid-based methods to solve the above elliptic interface problems accurately and robustly; however, we do not intend to have an exhaustive review here.

Besides the grid-based methods described above, the scientific computing community has shown an increased interest in solving elliptic interface problems using shallow or deep neural networks. Notice that the neural network approach for solving the interface problems has one apparent advantage over the grid-based methods; namely, it is completely mesh-free and can easily handle problems with complex interfaces or irregular domains. One obstacle for the neural network approach is that most of the network has a smooth activation function, so the resulting network is inherently smooth and is not a suitable ansatz for the interface problem. We list some related works in literature as follows. A deep Nitsche-type method~\cite{LM21} to solve elliptic interface problems with high-contrast discontinuous coefficients was developed in \cite{WZ20}. To deal with inhomogeneous boundary conditions, a shallow neural network to approximate the boundary conditions must be employed in advance. In \cite{GY22}, the authors proposed a deep unfitted Nitsche method for solving elliptic interface problems with high contrasts in high dimensions. Unlike using a single network, Wu and Lu~\cite{WL22} proposed an interfaced neural network that  decomposes the computational domain into two subdomains (one interface case), and each network is responsible for the solution on each subdomain. Then an extended multiple-gradient descent method was introduced to train the network. A similar piecewise deep neural network for elliptic interface problems was also introduced earlier in \cite{HHM22}. In the above neural network approaches, the network architectures usually have deep structures. Recently, the authors have proposed a discontinuity capturing shallow neural network (DCSNN)~\cite{HLL22} for solving elliptic interface problems with discontinuous solutions. By augmenting a coordinate variable to label different pieces of each subdomain, the DCSNN can be trained in a single physics-informed neural network~(PINN) framework~\cite{RPK19}. Meanwhile, we also used the idea proposed by E and Yu~\cite{EY18} and developed a completely shallow Ritz network for solving the elliptic interface problems by augmenting the level set function as an extra feature input in \cite{LCLHL22}. We found that it significantly improves the training effectiveness and accuracy. Notice that the major difference between DCSNN~\cite{HLL22} and the shallow Ritz network~\cite{LCLHL22} is that the former inherently represents a discontinuous function while the latter represents a continuous one.

In this paper, we propose a cusp-capturing physics-informed neural network for solving discontinuous-coefficient elliptic interface problems. The specific aim of this study is to introduce a network that can present continuous solutions, but with discontinuous first derivatives on interfaces. The smooth level set function augmented input in \cite{LCLHL22} cannot capture the derivative discontinuity sharply; thus, we augment a cusp-enforced level set function input to the network instead. Notice that, this new modified level set function does not change the interface position (i.e., zero level set). The rest of the paper is organized as follows. We present the formulation of the discontinuous-coefficient elliptic interface problems in Section~\ref{sec:ellipticFSI}. In Section~\ref{sec:cuspNN}, we propose a cusp-capturing neural network to solve the model problems. Numerical experiments are shown in Section~\ref{sec:numerics} to demonstrate the effectiveness of the proposed cusp-capturing technique and the accuracy of the present network, followed by some concluding remarks in Section~\ref{sec:conclusion}.

\section{Discontinuous-coefficient elliptic interface problems}
\label{sec:ellipticFSI}

We consider a $d$-dimensional discontinuous-coefficient second-order elliptic interface problem~\cite{BE87}. Let $\Omega\subset\mathbb{R}^d$ be a bounded domain and $\Gamma$ be an embedded $(d-1)$-dimensional $C^1$-interface separating $\Omega$ into two subdomains, $\Omega^-$ and $\Omega^+$, so $\Omega = \Omega^- \cup \Omega^+\cup \Gamma$. The equations of the problem subjected to the interface and boundary conditions are given as follows:
\begin{eqnarray}
\nabla\cdot\left(\beta(\bx)\nabla u(\bx)\right) - \alpha(\bx) u(\bx) &=&  f(\bx), \quad \bx \in \Omega^-\cup\Omega^+, \label{eq:varelliptic} \\
\dbblk{u}(\bx_\Gamma) \,=\, 0, \quad \dbblk{\beta\partial_n u}(\bx_\Gamma) &=& \rho(\bx_\Gamma), \quad \bx_\Gamma \in \Gamma, \label{eq:jump_cond}  \\  %
u(\bx_B) &=& g(\bx_B), \quad \bx_B \in\partial\Omega,\label{eq:bdc}
\end{eqnarray}
where $u(\bx)$ is the function to be solved,  $\rho(\bx_\Gamma)$ and $g(\bx_B)$ are given smooth functions, $\alpha(\bx)\ge 0$, $f(\bx)$ and $\beta(\bx)>0$ are also given but defined in a piecewise smooth manner across the interface $\Gamma$. We use $\partial_n u$ to denote the shorthand of normal derivative $\grad u\cdot\bn$, where $\bn$ is the unit normal vector pointing from $\Omega^-$ to $\Omega^+$  along the interface $\Gamma$. The notation $\dbblk{\cdot}$ represents the jump of a quantity across the interface (the one-sided limiting value approaching from $\Omega^+$ minus the one from $\Omega^-$).  For example,
\begin{equation}
\dbblk{\beta}(\bx_\Gamma) = \lim_{\bx\in\Omega^+,\,\bx\rightarrow\bx_\Gamma} \beta(\bx) - \lim_{\bx\in\Omega^-,\,\bx\rightarrow\bx_\Gamma}\beta(\bx) = \beta^+(\bx_\Gamma) - \beta^-(\bx_\Gamma), \label{eq:jump_def}
\end{equation}
where the superscripts ``$\pm$"  represent the limits of the function value on the interface.  Under this notation, the second interface condition in Eq.~(\ref{eq:jump_cond}) can be written explicitly as
\beqs
\dbblk{\beta\partial_n u}(\bx_\Gamma) & = & \beta^+(\bx_\Gamma)\partial_n u^+(\bx_\Gamma) - \beta^-(\bx_\Gamma)\partial_n u^-(\bx_\Gamma) \nonumber \\
& = & \beta^+(\bx_\Gamma)\partial_n u^+(\bx_\Gamma) - \beta^-(\bx_\Gamma)\partial_n u^+(\bx_\Gamma) + \beta^-(\bx_\Gamma)\partial_n u^+(\bx_\Gamma) - \beta^-(\bx_\Gamma)\partial_n u^-(\bx_\Gamma) \nonumber \\
&=& \dbblk{\beta}(\bx_\Gamma)\partial_n u^+(\bx_\Gamma)+ \beta^-(\bx_\Gamma)\dbblk{\partial_n u}(\bx_\Gamma) =\rho(\bx_\Gamma). \label{eq:jump_cond1}
\eeqs
One can immediately see that even with the case of $\dbblk{\beta}(\bx_\Gamma) = 0$,  the solution $u$ always has the property of $\dbblk{\partial_n u}(\bx_\Gamma)\neq 0$ as long as $\rho(\bx_\Gamma)\neq 0$. Along with the first interface condition
$\dbblk{ u}(\bx_\Gamma) = 0$ in Eq.~(\ref{eq:jump_cond}), we can conclude that the solution $u$ is continuous over the domain $\Omega$ but its normal derivative has jump discontinuity across the interface $\Gamma$.

We would also like to point out that although here we focus only on the Dirichlet-type boundary condition~(\ref{eq:bdc}), one can apply the present method to the Neumann or Robin-type boundary condition with no difficulty. In this paper, we aim to find the solution to Eqs.~(\ref{eq:varelliptic})-(\ref{eq:bdc}) using machine learning techniques in the spirit of physics-informed neural networks~\cite{RPK19}, as introduced in the next section.

\section{A cusp-capturing physics-informed neural network}
\label{sec:cuspNN}

As mentioned before, the solution of Eqs.~(\ref{eq:varelliptic})-(\ref{eq:bdc}) is continuous in the domain $\Omega$ but has a jump discontinuity to its normal derivative on the interface $\Gamma$. The universal approximation theorems~\cite{Cybenko89,Hornik89,Pinkus99} guarantee the applicability of approximating such continuous solutions using artificial neural networks. However, a neural network with differentiable activation functions is undoubtedly smooth, thus it is unlikely to capture the present solution with cusps (the partial derivatives are not continuous) in an accurate manner. More precisely, locating and fitting derivative discontinuities in neural network solutions is challenging. Since the partial derivative jumps occur at the interface, it is natural to include the interface position as a feature input in the network architecture. In \cite{LCLHL22}, we proposed a shallow Ritz-type method to solve similar interface problems (taking $\beta=1$) as Eqs.~(\ref{eq:varelliptic})-(\ref{eq:bdc}) in which we add the level set function of the interface as a feature input to the network. That is, we
use a neural network of the form $U(\bx, z=\phi(\bx))$ to approximate the solution $u(\bx)$ of the problem, where $\phi(\bx)$ is the level set function defined in the whole domain $\Omega$. Here, the interior and exterior region are defined as
$\Omega^- = \{\bx \in\mathbb{R}^d\, |\phi(\bx)<0\}$ and $\Omega^+ = \{\bx \in\mathbb{R}^d\, |\phi(\bx)>0\}$, respectively, and the zero level set gives the position of the interface $\Gamma$, i.e., $\Gamma = \{\bx \in\mathbb{R}^d\,|\, \phi(\bx)=0\}$. With this level set function augmentation, we found that it significantly improves the training effectiveness and accuracy. However, since the level set function is smooth, and the neural network function $U$ is smooth due to the use of a smooth activation function, the resulting neural network solution $u(\bx)=U(\bx, z)=U(\bx, \phi(\bx))$ remains smooth. That is, the gradient of $u$
\beq
\grad u = \grad_\bx U + \pd _z U \, \grad \phi,
\label{eq:cusp0}
\eeq
is continuous so the normal derivative jump $\dbblk{\partial_n u}=0$ across the interface $\Gamma$.
Here, $\nabla_\bx U\in\mathbb{R}^d$ represents a vector with partial derivatives of $U$ with respect to the components in $\bx$, and $\partial_zU$ is the partial derivative of $U$ with respect to $z$. We also suppress the notation of $\bx$ in the gradients of $u$ and $\phi$ since they both are functions of $\bx$. Thus, if we want to require $\grad u$ to be discontinuous across the interface then $\grad \phi$ should be discontinuous too. Therefore, we need to modify the original smooth level set function accordingly.

\subsection{Cusp-enforced level set function augmentation}
As mentioned above, we need to modify the level set function so that its gradient is discontinuous across the interface without changing the zero level set. This can be done easily by taking the absolute value of the level set function; that is, we define $\phi_a(\bx)= |\phi(\bx)|$. We therefore call this $\phi_a$ as a cusp-enforced level set function since it is non-differentiable at the interface $\Gamma$. Furthermore, one can immediately derive that this cusp-enforced level set function has the gradient jump as $\dbblk{\grad \phi_a}(\bx_\Gamma) = 2 \grad \phi(\bx_\Gamma), \bx_\Gamma \in \Gamma$. Note that, the above jump condition is evaluated by the limiting values from both sides of the interface where $\nabla\phi_a$ is well-defined.  With this modified level set function,  we now define a new neural network solution in the form as $u(\bx) =U(\bx, z)=U(\bx, \phi_a(\bx))$. Since the neural network function $U$ is smooth, calculating the derivatives of the network $U$ with respect to its input variables $\bx$ and $z$ via automatic differentiation~\cite{GW08} has no problem at all. Thus, the gradient jump of $u$ across the interface can be computed directly from Eq.~(\ref{eq:cusp0}) as
\beq
\dbblk{\grad u}(\bx_\Gamma) = \pd _z U \dbblk{\grad \phi_a}(\bx_\Gamma) = 2 \pd _z U  \grad \phi(\bx_\Gamma).
\eeq
Notice that, in the above implementation we have used $\dbblk{\grad_\bx U}(\bx_\Gamma)=0$ since $U$ is smooth.
By multiplying the normal vector $\bn=\grad \phi/ \|\grad \phi \|$ to the above equation, we obtain the following normal derivative jump of $u$ as
\beq
 \dbblk{\pd_n u}(\bx_\Gamma) = 2 \pd _z U \|\grad \phi(\bx_\Gamma) \|.
\eeq
Therefore, the neural network solution $U$ is capable of capturing the cusp behavior of the solution in Eqs.~(\ref{eq:varelliptic})-(\ref{eq:bdc}) even if the network function $U(\bx, z)$ is smooth across its entire $\mathbb{R}^{d+1}$ domain.

By using the relation $\grad u = \grad_\bx U + \pd _z U \, \grad \phi_a$ in $\Omega^\pm$, one can explicitly write the following equation after careful calculations
\begin{eqnarray}
\label{eq:lapU}
\nabla\cdot\left(\beta \nabla u\right)
&=& \beta\left(\Delta_\bx U + 2\grad\phi_a\cdot\grad_\bx\left( \pd_{z}U\right)+\|\grad\phi\|^2\pd_{zz}U+\pd_z U\Delta\phi_a\right) \\
&& \, + \quad\grad\beta\cdot\left(\grad_\bx U + \pd_{z}U \grad\phi_a \right), \nonumber
\end{eqnarray}
where $\Delta_\bx$ is the Laplace operator concerning only the variable $\bx$.

Now, Eqs.~(\ref{eq:varelliptic})-(\ref{eq:bdc}) can be rewritten in terms of $U$ as follows. For succinctness, we introduce the notation $\mathcal{L}_{\beta,\phi_a}U$ to represent the right-hand side of Eq.~(\ref{eq:lapU}) so that Eq.~(\ref{eq:varelliptic}) is rewritten to the following
\begin{equation}
\label{eq:poissonU}
\mathcal{L}_{\beta,\phi_a}U(\bx,\phi_a(\bx)) - \alpha(\bx)U(\bx,\phi_a(\bx)) = f(\bx), \quad \bx \in \Omega^+\cup\Omega^-.
\end{equation}
Using the fact that $\dbblk{\partial_n\phi_a} (\bx_\Gamma)= 2 \grad \phi(\bx_\Gamma)\cdot \bn = 2 \| \grad \phi(\bx_\Gamma) \|$,  we can also rewrite the interface condition $\dbblk{\beta\partial_n u}(\bx_\Gamma) =\rho(\bx_\Gamma)$  in Eq.~(\ref{eq:jump_cond1}) as
\begin{equation}
\label{eq:jump_condU}
 \dbblk{\beta}(\bx_\Gamma)\partial_n U + (\beta^+(\bx_\Gamma) +\beta^-(\bx_\Gamma)) \partial_{z}U\, \| \grad \phi(\bx_\Gamma) \| = \rho(\bx_\Gamma) \quad \bx_\Gamma \in \Gamma,
\end{equation}
where $\partial_n U = \grad_\bx U\cdot\bn$. Notice that $\dbblk{u}(\bx_\Gamma)=0$ is automatically satisfied  since $U$ is a continuous function. The associated boundary condition~(\ref{eq:bdc}) reads
\begin{equation}
\label{eq:bdc1}
U(\mathbf{x}_B,\phi_a(\mathbf{x}_B)) = g(\mathbf{x}_B) \quad \bx_B \in\partial\Omega.
\end{equation}

The remaining task is to train the network to simultaneously satisfy Eq.~(\ref{eq:poissonU}), the jump condition (\ref{eq:jump_condU}), and the boundary condition (\ref{eq:bdc1}) with appropriate loss function.

\subsection{Physics-informed neural networks}
\begin{figure}[h]
\begin{center}
\includegraphics[scale=0.5]{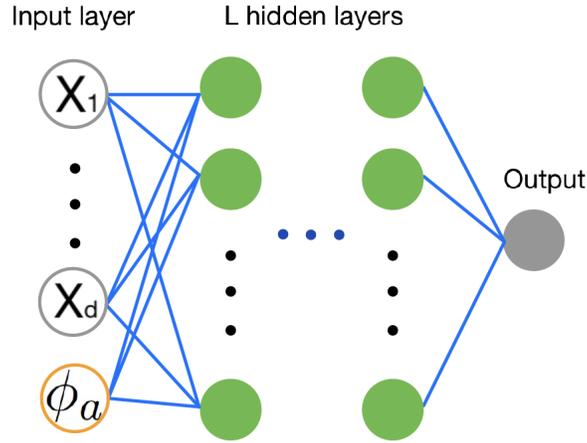}
\caption{Diagram of the $L$-hidden-layer network structure.}\label{fig:shallow&deepNN}
\end{center}
\end{figure}

In this subsection, we present a physics-informed neural network to approximate the solution $U(\bx, \phi_a(\bx))$ for Eqs~.(\ref{eq:poissonU})-(\ref{eq:bdc1}). The convergence of PINNs for linear elliptic PDEs was studied recently in ~\cite{SDK20}. Figure~\ref{fig:shallow&deepNN} presents the structure of a $L$-hidden-layer feed-forward fully connected neural network where $(\bx,\phi_a(\bx))^T\in\mathbb{R}^{d+1}$ represents the $d+1$ feature input of the network (recall that $ \phi_a(\bx)$ is the cusp-enforced level set function). We label the input layer as layer $0$ and denote the feature input as $\bv^{[0]}=(\bx,\phi_a(\bx))^T$. The output at the $\ell$-th hidden layer with $N_\ell$ neurons, denoted as $\bv^{[\ell]}\in\mathbb{R}^{N_\ell}$, presents an affine mapping of the output of layer $\ell-1$ (i.e., $\bv^{[\ell-1]}$) followed by an action of the activation function $\sigma$ in a componentwise manner as
\begin{equation}
\bv^{[\ell]} = \sigma\left(W^{[\ell]}\bv^{[\ell-1]}+\bb^{[\ell]}\right), \quad \ell=1,\cdots,L, \label{eq:DNN}
\end{equation}
where the matrix $W^{[\ell]}\in\mathbb{R}^{N_\ell\times N_{\ell-1}}$ contains the weights connecting the structure from layer $\ell-1$ to layer $\ell$, and $\bb^{[\ell]}\in\mathbb{R}^{N_\ell}$ is the bias vector at layer $\ell$. Finally, we denote the output of this multiple-hidden-layer network as
\begin{equation}
U_{\mathcal{N}}(\bx, \phi_a(\bx);\bft) = W^{[L+1]}\mathbf{v}^{[L]},
\end{equation}
where $W^{[L+1]}\in\mathbb{R}^{1\times N_L}$. The notation $\bft$ denotes the vector collecting all trainable parameters (including all the weights and biases) so the dimension of $\bft$ is the total number of parameters in the network that can be easily counted as $N_\bft = N_{L}+\Sigma^{L}_{\ell=1}(N_{\ell-1}+1)N_\ell$.

In the training process, we select $M_I$ points in the region of $\Omega^-\cup\Omega^+$, $\left\{\bx^i\right\}_{i=1}^{M_I}$, $M_\Gamma$ points on the interface $\Gamma$, $\left\{\bx^i_\Gamma\right\}_{i=1}^{M_\Gamma}$, and $M_B$ points on the domain boundary $\partial\Omega$, $\left\{\bx^i_B\right\}_{i=1}^{M_B}$, so totally $M=M_I+M_\Gamma+M_B$ training points. Under the physics-informed framework, we hereby define the loss function as the mean squared error of the residual of differential equation~(\ref{eq:poissonU}), the jump condition~(\ref{eq:jump_condU}), and the boundary condition (\ref{eq:bdc1}) as
\begin{eqnarray}
\label{eq:loss}
\Loss(\bft) &=& \frac{1}{M_{I}}\sum_{i=1}^{M_{I}} \abs{L_I(\bx^i,\phi_a(\bx^i);\bft)}^{2} + \frac{c_\Gamma}{M_\Gamma}\sum_{i=1}^{M_\Gamma} \abs{ L_\Gamma(\bx^i_\Gamma,0;\bft)}^{2} \nonumber\\ && \quad + \frac{c_B}{M_{B}}\sum_{i=1}^{M_{B}} \abs{L_B(\bx^i_B,\phi_a(\bx^i_B);\bft)}^{2}.
\end{eqnarray}
where the residual error $L_I$, interface condition error $L_\Gamma$, and boundary condition error $L_B$, are shown respectively as follows:
\begin{eqnarray}
L_I(\bx,\phi_a(\bx);\bft) &=& \mathcal{L}_{\beta,\phi_a}U_\mathcal{N}(\bx,\phi_a(\bx);\bft) - \alpha(\bx)U_\mathcal{N}(\bx,\phi_a(\bx);\bft) - f(\bx), \label{eq:lossi}\\
L_\Gamma(\bx_\Gamma,0;\bft) &=& \dbblk{\beta}(\bx_\Gamma)\partial_n  U_\mathcal{N}(\bx_\Gamma,0;\bft)\, + (\beta^+(\bx_\Gamma) +\beta^-(\bx_\Gamma)) \partial_{z}U_\mathcal{N}(\bx_\Gamma,0;\bft)\, \| \grad \phi(\bx_\Gamma) \| \nonumber \\
& &- \rho(\bx_\Gamma), \label{eq:lossg} \\
L_B(\bx_B,\phi_a(\bx_B);\bft) &=& U_\mathcal{N}(\bx_B,\phi_a(\bx_B);\bft)-g(\bx_B). \label{eq:lossb}
\end{eqnarray}
The constants $c_\Gamma$ and $c_B$ appeared in the loss function~(\ref{eq:loss}) are chosen to balance the contribution of the terms related to the interface jump condition~(\ref{eq:jump_condU}) and boundary condition~(\ref{eq:bdc1}), respectively. In latter numerical experiments,  we might need to use network with smooth level set function $\phi$ augmentation $U_\mathcal{N}(\bx,\phi(\bx);\bft)$ for comparison purpose. In that case,
the interface error loss in Eq.~(\ref{eq:lossg}) should be replaced (can be easily derived) by
\begin{equation}
L_\Gamma(\bx_\Gamma,0;\theta) = \dbblk{\beta}(\bx_\Gamma)\left(\partial_n  U_\mathcal{N}(\bx_\Gamma,0;\bft)+ \partial_{z}U_\mathcal{N}(\bx_\Gamma,0;\bft)\| \grad \phi(\bx_\Gamma) \|\right)- \rho(\bx_\Gamma). \label{eq:lossg2}
\end{equation}
Meanwhile, throughout the rest of paper, we use the Levenberg-Marquardt (LM) algorithm~\cite{More78} as the optimizer to train the network, and  use the notation $u_\mathcal{N}$ to denote the network prediction solution.
\vspace{3mm}\\
{\bf Remark.} The cusp-capturing PINN is designed for solving elliptic interface problems where the solution is continuous but the derivatives have jumps. The present method can be easily extended to handle problems with non-zero solution jumps. If the solution is discontinuous across the interface, we can incorporate an additional supervised learning task for solution jump approximation and the remaining part of the solution can be found by the cusp-capturing PINN. To see this, suppose we want to solve Eqs.~(\ref{eq:varelliptic})-(\ref{eq:bdc}) but  with nonzero solution jump $\dbblk{u}(\bx_\Gamma)=\lambda(\bx_\Gamma), \forall \bx_\Gamma \in \Gamma$ instead. We first write the solution  as $u(\bx) = v(\bx) + w(\bx)$ in which we assume $v(\bx)$ has the jump  discontinuity $\dbblk{v}(\bx_\Gamma)=\lambda(\bx_\Gamma)$ so $w(\bx)$ is continuous ($\dbblk{w}(\bx_\Gamma)=0$).  We further assume $v(\bx)$ has the form
\begin{align}\label{Eq:v_V}
v(\bx) =
\left\{
\begin{array}{ll}
V(\bx)  & \bx \in \Omega^-,\\
0          & \bx \in \Omega^+,
\end{array}\right.
\end{align}
so the jump $\dbblk{v}(\bx_\Gamma)=-V(\bx_\Gamma)=\lambda(\bx_\Gamma)$ for $\bx_\Gamma \in \Gamma$. The construction of $V(\bx)$ will become clear later. Substituting the expression of $u(\bx)$ into
Eqs.~(\ref{eq:varelliptic})-(\ref{eq:bdc}), one can immediately obtain the equations for $w(\bx)$ as
\begin{eqnarray}
\nabla\cdot\left(\beta(\bx)\nabla w(\bx)\right) - \alpha(\bx) w(\bx) &=& \left\{\bary{ll} f(\bx)-\nabla\cdot\left(\beta(\bx)\nabla V(\bx)\right)+\alpha(\bx)V(\bx), & \bx \in \Omega^-, \\ f(\bx), & \bx\in\Omega^+ \eary\right.
\label{eq:varelliptic2} \\
\dbblk{w}(\bx_\Gamma) \,=\, 0, \quad \dbblk{\beta\partial_n w}(\bx_\Gamma) &=& \rho(\bx_\Gamma)+ \beta^-(\bx_\Gamma)\partial_n V(\bx_\Gamma), \quad \bx_\Gamma \in \Gamma, \label{eq:jump_cond2}  \\  %
w(\bx_B) &=& g(\bx_B), \quad \bx_B \in\partial\Omega.\label{eq:bdc2}
\end{eqnarray}
Note that, the flux jump in Eq.~(\ref{eq:jump_cond2}) is obtained by the fact $\dbblk{\beta\partial_n v}(\bx_\Gamma)=-\beta^-(\bx_\Gamma)\partial_n V(\bx_\Gamma)$. The above equations (\ref{eq:varelliptic2})-(\ref{eq:bdc2}) can be solved by the present cusp-capturing PINN since the solution $w(\bx)$ now is continuous.

The remaining question is how to construct the function $V(\bx)$ so that $V(\bx_\Gamma)=-\lambda(\bx_\Gamma)$ for $\bx_\Gamma \in \Gamma$.  Here, we simply adopt a shallow (one-hidden-layer) fully-connected feedforward neural network to approximate $V$  by supervised learning. That is, we randomly choose $M_\Gamma$ points $\left\{\bx^i_\Gamma\right\}_{i=1}^{M_\Gamma}$ on the interface $\Gamma$, and minimize the corresponding mean squared error loss as
\beq
	\Loss(\tilde{\theta})=  \frac{1}{M_\Gamma} \sum_{i=1}^{M_\Gamma}\left(V(\bx^i_\Gamma;\tilde{\theta})+\lambda(\bx^i_\Gamma)\right)^2, 
\eeq		
where $\tilde{\theta}$ denotes the vector collecting the trainable weights and biases used in the network.

\section{Numerical results}
\label{sec:numerics}

In this section, we aim to demonstrate the capability of the present neural network method for solving elliptic interface problems, Eqs.~(\ref{eq:varelliptic})-(\ref{eq:bdc}). We set the penalty constants in the loss function $c_B=c_\Gamma=1$ to focus on the accuracy check of the present cusp-capturing technique. The merit of the proposed cusp-capturing PINN is to allow one to use a smooth neural network $U(\bx, z)$ to learn the non-smooth solution, $u(\bx)$, through the relation $u(\bx) = U(\bx,z=\phi_a(\bx))$. The only requirement of the choice of activation function is subject to the $C^2$-regularity of $u(\bx)$ in each subdomain. Thus, we simply choose the sigmoid function, $\sigma(x)=\frac{1}{1+e^{-x}}$ as our activation function. For the following numerical examples, we employ different depth networks (from $1$ to $L$  hidden layers) with equal number of neurons in each hidden layer $N_1 = N_2 = \cdots = N_L = N$. The training and test data points are generated by the Latin hypercube sampling algorithm~\cite{MBC79}, which effectively avoids the clustering of data points at some specific locations so resulting in a nearly random sampling. To measure the accuracy of the network solution, we choose $M_{test}$ points (different from the training points) in $\Omega$ to calculate the relative $L^\infty$ and $L^2$ errors defined respectively as $\|u_\mathcal{N}-u\|_\infty/\|u\|_\infty$ and $\|u_\mathcal{N}-u\|_2/\|u\|_2$, where
\begin{align*}
\|u\|_\infty = \max_{1\leq i \leq M_{test}}|u(\bx^i)|, \quad
\|u\|_2 = \sqrt{\frac{1}{M_{test}}\sum_{i=1}^{M_{test}}(u(\bx^i))^2}.
\end{align*}
In general, we set $M_{test}=100 M$, where $M$ is the total number of training points. Since the predicted results will vary slightly for each experiment (it is affected by the randomness of the training and test data points, and the initialization of trainable parameters), we show the average value of the errors and losses over 5 trial runs.

In the training procedure, we use the Levenberg-Marquardt (LM) algorithm  as our optimizer and update the damping parameter $\mu$ by the strategies introduced in~\cite{TS12}. The training is stopped when the loss value $\Loss(\bft)$ is below a threshold $\epsilon_\theta$ (problem dependent) or the maximum iteration (training) step $epoch=3000$ is reached. All trials are run on a desktop equipped with one NVIDIA GeForce RTX3060 GPU. We implement the cusp-capturing PINN architecture using Pytorch (v1.13)~\cite{PGMetel19} and all trainable parameters (weights and biases) are initialized using Pytorch default settings. The source codes used throughout this paper are available on GitHub at https://github.com/teshenglin/cusp\_capturing\_PINN.

\paragraph{Example 1.} As the first example, we demonstrate the cusp-capturing capability for the present network by considering the following one-dimensional Poisson equation on an interval $\Omega=[0,1]$ with an interface point at $x_\Gamma=\frac{1}{3}$:
\beqs
\frac{d^2 u}{d x^2} & = & 0, \quad x \in (0,1)\backslash \{x_\Gamma\}, \label{1D_Poisson}\\
\dbblk{u}(x_\Gamma) =0, \quad \dbblk{\frac{du}{dx}}(x_\Gamma) & = & 1, \\
u(0)=u(1)& = & 0.
\eeqs
The exact solution of the above problem can be easily derived as
\begin{equation}
u(x) = \left\{ \begin{array}{ll}
(x_\Gamma-1)x, &  x\in [0,x_\Gamma), \\
x_\Gamma(x-1), &  x\in [x_\Gamma,\,1],
\end{array}\right.
\end{equation}
where the cusp appears exactly at the interface $x_\Gamma$. We thus choose $\phi(x)=x-x_\Gamma$ as the smooth level set function so that $\phi_a(x) = \abs{x-x_\Gamma}$ represents the cusp-enforced level set function.

For the neural network in this test, we use a completely shallow network structure ($L=1$) with $N$ neurons in the hidden layer; here, the input dimension is two, one for $x$ and the other for the augmented feature input $\phi_a$. The number of overall training data points is $M=M_I+3$, including $M_I$ points in the interval $(0,1)$, two points ($M_B=2$) at the  boundary, and one point ($M_\Gamma=1$) at the interface. We use only $2$ neurons in the hidden layer and $13$ training points, that is, $(N, M) = (2, 13)$. After completing the training process, we use $M_{test}=1000$ test points to examine the predicted accuracy of the network solution.

\begin{figure}[!]
\begin{center}
\includegraphics[scale=0.34]{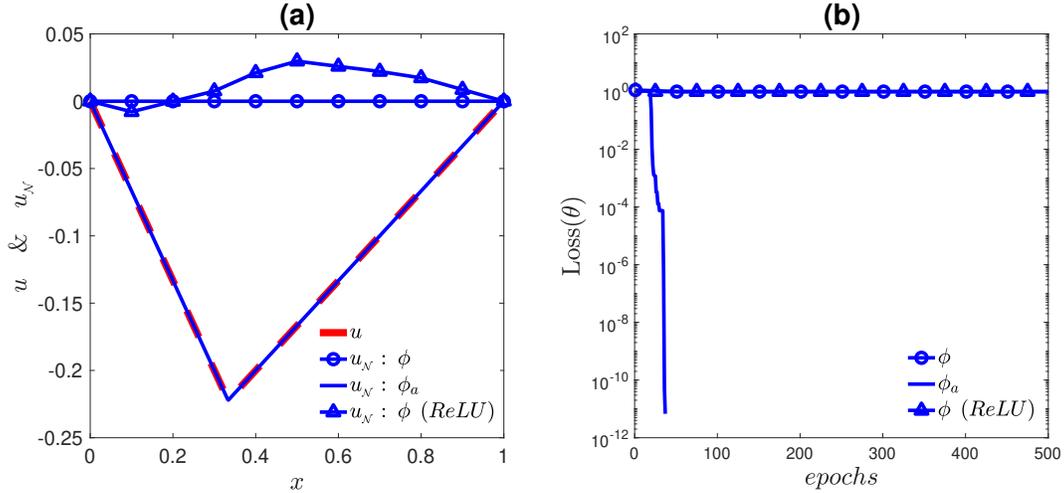}
\caption{(a) The profiles of the exact solution $u$, the network solutions $u_\mathcal{N}$ with augmented input $\phi$ and $\phi_a$, and the network solution using ReLU activation function with $\phi$ augmented input. (b) The corresponding losses in (a).}
\label{fig:ex01_Poisson1D_greenfun_ctsphi&absphi&ctsphiReLU}
\end{center}
\end{figure}

Figure~\ref{fig:ex01_Poisson1D_greenfun_ctsphi&absphi&ctsphiReLU}(a) shows the profiles of the exact solution $u$ (denoted by the red-dashed line) and the network-predicted solution $u_\mathcal{N}$ with augmented input $\phi_a$ (solid line). One can immediately see that the $\phi_a$ input network solution captures the cusp sharply where the $L^\infty$ error achieves $\|u-u_\mathcal{N}\|_\infty=7.01\times 10^{-8}$. Meanwhile, the corresponding loss drops significantly within just $40$ epochs, as shown in panel (b) of the figure.

Then we test to see if the solution can be learned by using a level set function augmented input (not the cusp-enforced one); that is, we assume $u_\mathcal{N} = U_\mathcal{N}(\bx, \phi(\bx))$. 
 We train the network with $(N, M) = (20, 103)$. The learned solution is shown in Figure~\ref{fig:ex01_Poisson1D_greenfun_ctsphi&absphi&ctsphiReLU}(a) (denoted by ``$\circ$") and the corresponding loss is presented in (b). It turns out that the $\phi$ input network learns a completely wrong solution $u_\mathcal{N} \approx 0$. This result is not surprising, since this network solution is inherently smooth, so all the jumps are zero, which gives $L_\Gamma(\bx_\Gamma,0;\theta) = - \rho(\bx_\Gamma)$ that is independent of the trainable parameters $\theta$. So this smooth neural network tries to minimize only the residual error and boundary error, that is, to learn a solution with zero second-order derivative and zero boundary condition. The loss for this $\phi$ input network shown in panel (b) is dominated by the interface loss $L_\Gamma$ that gives an $O(1)$ value throughout the whole training process.

Meanwhile, one may wonder if a feed-forward network using the ReLU activation function with augmented smooth level set function $\phi$ can work due to the cusp-like profile of the ReLU function. Notice that, the ReLU function is linear so a shallow network (one hidden layer) with ReLU activation can learn the differential equation (\ref{1D_Poisson}) with zero loss (i.e $L_I(x, \phi(x), \theta)=0$). However, it seems to be difficult to locate the cusp singularity for such a network which we can see from the solution profile (denote by ``$\triangle$")  in Figure~\ref{fig:ex01_Poisson1D_greenfun_ctsphi&absphi&ctsphiReLU}(a). Again, like the sigmoid activation function with $\phi$ augmented input, the corresponding loss (also see in Figure~\ref{fig:ex01_Poisson1D_greenfun_ctsphi&absphi&ctsphiReLU}(b)) remains to be $O(1)$ which leads to unsuccessful training. As discussed in \cite{WL22}, a single network with non-differentiable activation usually does not satisfy the differential requirement in high-dimensional interface problems. As a result, the cusp singularity obtained by the network does not coincide with the given interface. This is exactly what we see from Figure~\ref{fig:ex01_Poisson1D_greenfun_ctsphi&absphi&ctsphiReLU}(a) even in a one-dimensional case.

\paragraph{Example 2.}
As the second example, we consider an elliptic equation with a piecewise-constant coefficient defined in the two-dimensional domain $\Omega=[-1,1]\times[-1,1]$. The embedded interface $\Gamma$ is described by the zero level set of the function $\phi(x,y) = \frac{x^2}{0.5^2} + \frac{y^2}{0.5^2}-1$, separating $\Omega$ into the inner ($\Omega^-$) and outer ($\Omega^+$) regions. We choose the exact solution $u$ and the coefficient $\beta$, respectively, as
\begin{equation}
u(x,y) = \left\{ \begin{array}{ll}  1-\exp{\left(\frac{1}{\eta}\left(\frac{x^2}{0.5^2}+\frac{y^2}{0.5^2}-1\right)\right)}, &  (x,y)\in\Omega^-,  \\ -\gamma \ln\left(\frac{x^2}{0.5^2} + \frac{y^2}{0.5^2}\right), &  (x,y)\in\Omega^+, \end{array}\right.
\end{equation}
and
\begin{equation*}
\beta(x,y)=\left\{ \begin{array}{ll}  \beta^-, & \quad(x,y)\in \Omega^-, \\ \beta^+, & \quad(x,y)\in \Omega^+,\end{array}\right.	
\end{equation*}
where the parameter $\eta=\beta^-/\beta^+$ represents the ratio of $\beta^-$ to $\beta^+$.  (Here, we fix $\beta^+=1$ and adjust $\eta$ to control the contrast of the coefficients.)  One can immediately see that the solution is continuous across the interface $\Gamma$ but its normal derivative has jump discontinuity as $\dbblk{\beta\partial_n u}=-4(\gamma-1)$. The corresponding right-hand side function $f$ can be calculated directly from Eq.~(\ref{eq:varelliptic}) and the boundary condition $g$ is given by the exact solution $u$ on $\partial\Omega$. We introduce a number $M_0$ which can be regarded as the grid number used in each spatial dimension as in traditional grid-based methods so the training data set includes $M_I = M^2_0$ points in $\Omega^-\cup\Omega^+$, $M_\Gamma = 3M_0$ points on the interface $\Gamma$, and $M_B = 4M_0$ points on the boundary $\partial\Omega$, respectively. Thus, the total training points $M = M_0^2 + 3M_0 + 4M_0$. 

Next, we will discuss some numerical issues about the implementation of cusp-capturing strategy, including the accuracy study of shallow neural networks with different number of neurons and training points, and the comparisons of different optimizers and different augmented inputs.

\begin{figure}[!]
\begin{center}
\includegraphics[width=\textwidth]{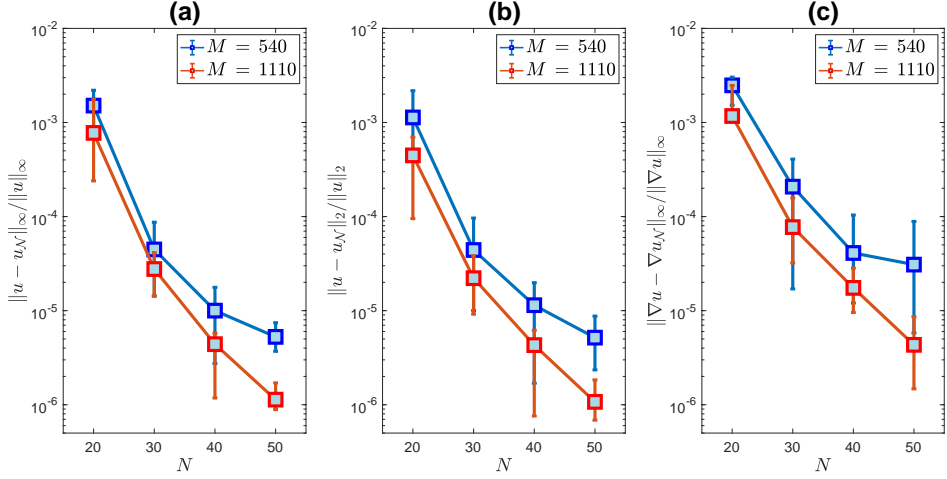}
\caption{Error bar plots associated with Table 1. Each bar represents the errors over 5 trial runs. (a) Relative $L^\infty$ error of $u_\mathcal{N}$; (b) Relative $L^2$ error of $u_\mathcal{N}$; (c) Relative $L^\infty$ error of $\nabla u_\mathcal{N}$.}
\label{fig:ex02_tab01_errbar}
\end{center}
\end{figure}

\paragraph{\emph{Accuracy check: shallow neural networks with different number of neurons and training points.}} The first experiment aims to study the number of neurons and training points needed to get satisfactory results. To test whether the proposed method works for different types of boundary condition, we impose the Dirichlet boundary condition at $x=\pm1$ and the Neumann boundary condition at $y=\pm1$. We choose $\alpha=1$, $\eta=10$, $\gamma=2$, and fix $L=1$ such that the neural network is completely shallow.
\begin{table}[h!]
\begin{center}
   \begin{tabular}{c|c|ccccc}
   \hline
    $(M_0,M)$ & $(N, N_\bft)$ & $\frac{\|u_\mathcal{N}-u\|_\infty}{\|u\|_\infty}$ & $\frac{\|u_\mathcal{N}-u\|_2}{\|u\|_2}$ & $\frac{\|\nabla u_\mathcal{N}-\nabla u\|_\infty}{\|\nabla u\|_\infty}$ & $\Loss(\bft)$ & \\
    \hline
 	    	& $(20,100)$ & $1.13\times10^{-3}$ & $2.82\times10^{-3}$ & $2.48\times10^{-3}$ & $2.23\times10^{-4}$  & \\
 $(20,540)$ & $(30,150)$ & $4.41\times10^{-5}$ & $2.33\times10^{-4}$ & $2.08\times10^{-4}$ & $1.61\times10^{-7}$  & \\
 			& $(40,200)$ & $1.14\times10^{-5}$ & $4.56\times10^{-5}$ & $4.10\times10^{-5}$ & $2.29\times10^{-9}$  & \\
 			& $(50,250)$ & $5.17\times10^{-6}$ & $3.69\times10^{-5}$ & $3.09\times10^{-5}$ & $1.04\times10^{-10}$ & \\
   \hline
 			& $(20,100)$ & $7.75\times10^{-4}$ & $4.47\times10^{-4}$ & $1.17\times10^{-3}$ & $7.12\times10^{-5}$  & \\
$(30,1110)$ & $(30,150)$ & $2.77\times10^{-5}$ & $2.21\times10^{-5}$ & $7.73\times10^{-5}$ & $8.79\times10^{-8}$  & \\
 			& $(40,200)$ & $4.40\times10^{-6}$ & $4.30\times10^{-6}$ & $1.75\times10^{-5}$ & $2.97\times10^{-9}$  & \\
 			& $(50,250)$ & $1.13\times10^{-6}$ & $1.07\times10^{-6}$ & $4.33\times10^{-6}$ & $1.17\times10^{-10}$ & \\
   \hline
    \end{tabular}
\caption{Relative errors of $u$ and $\nabla u$, and training losses for the shallow network solution with different number of neurons $N$ and training points $M$. Here, $\alpha=1$, $\eta=10$, and $\gamma=2$ in Example~2.} \label{tab:ex02_eta10_N&Mcvg}
    \end{center}	
\end{table}
Table~\ref{tab:ex02_eta10_N&Mcvg} shows the relative $L^\infty$ and $L^2$ errors between the network solution $u_\mathcal{N}$ and the exact solution $u$ when using different numbers of neurons $N$ and training points $M$. Also, we examine the relative $L^\infty$ error of $\nabla u_\mathcal{N}$ by the formula $\|\nabla u_\mathcal{N}-\nabla u\|_\infty/\|\nabla u\|_\infty$ with the definition $\|\nabla{u}\|_\infty=\frac{1}{2}\left(\|\frac{\pd u}{\pd x}\|_\infty+\|\frac{\pd u}{\pd y}\|_\infty\right)$. Notice that since the network has only one hidden layer, the overall number of trainable parameters is $N_\theta = N (d+3)= 5N$ for this two-dimensional problem. The corresponding final loss values are also shown in the table. One can see that the present model can achieve a prediction accuracy of about $0.1\%$ in relative $L^\infty$ and $L^2$ errors even using one hidden layer with merely $N = 20$ neurons. As we increase the number of neurons from $N=20$ to $N=50$, the relative error decreases from the magnitude $O(10^{-3})$ to $O(10^{-6})$, and the loss drops from $O(10^{-4})$ to $O(10^{-10})$ accordingly. In addition, one can also see that all relative errors decrease by increasing the number $M_0=20$ to $M_0=30$ (same as increasing the number of total training points $M$). From this numerical experiment, we conclude that the solution errors can indeed be reduced by increasing the number of neurons or training points, which provides an informal evidence for the numerical convergence of the present method. The errors for the solution gradient show a similar convergence trend as the solution errors. In addition, since the derivatives are computed by automatic differentiation, the relative $L^\infty$ errors of the gradient seem to have almost the same order of magnitude as the ones of the solution itself. We also present the error bar plots of 5 trail runs associated with Table~\ref{tab:ex02_eta10_N&Mcvg} in Figure~\ref{fig:ex02_tab01_errbar}.

\begin{figure}[!]
\begin{center}
\includegraphics[width=\textwidth]{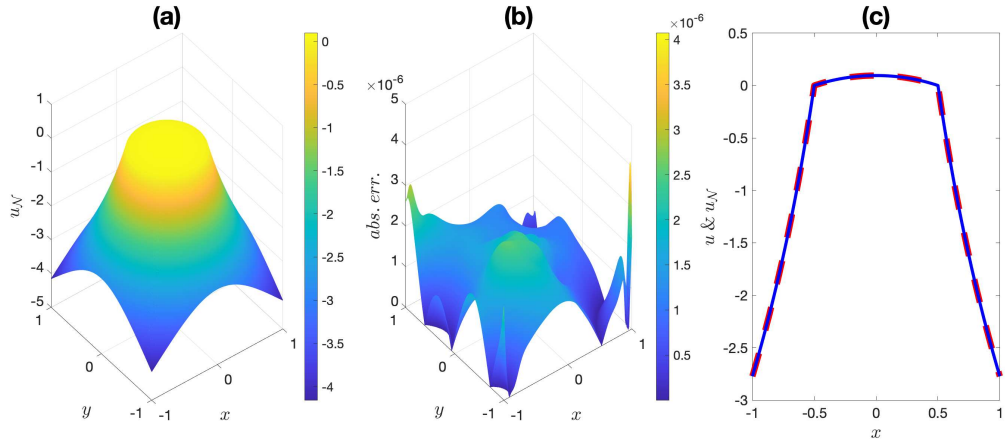}
\caption{ (a) The solution profile of $u_\mathcal{N}$; (b) Absolute error $\abs{u_\mathcal{N}-u}$; (c) Cross-sectional view of $u_\mathcal{N}$ (blue-solid line) and $u$ (red-dashed line) along the line $y=0$. The figure is the case when $(M_0,M)= (30, 1110)$ and $(N,N_\bft)=(50,250)$ in Table 1.}
\label{fig:ex02_eta10_uS_abserr_crossx_S40M1110}
\end{center}
\end{figure}

We depict the solution profile $u_\mathcal{N}$ in Figure~\ref{fig:ex02_eta10_uS_abserr_crossx_S40M1110}(a), the corresponding absolute error $\abs{u_\mathcal{N}-u}$ in Figure~\ref{fig:ex02_eta10_uS_abserr_crossx_S40M1110}(b), and the cross-sectional view of $u_\mathcal{N}$ and $u$ along the line $y=0$  in Figure~\ref{fig:ex02_eta10_uS_abserr_crossx_S40M1110}(c). One can clearly see that the cusps on the interface are accurately captured and the largest error occurs at the domain boundary rather than on the interface, which indicates the effectiveness of the present network model.

\begin{figure}[!]
\begin{center}
\includegraphics[width=\textwidth]{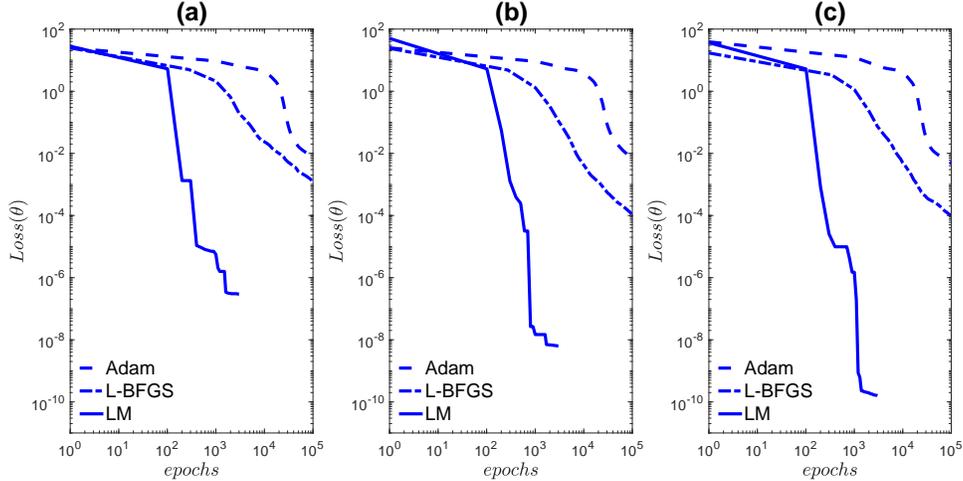}
\caption{Comparison of loss evolutions using different optimizers: Adam (dashed line), L-BFGS (dashed-dotted line), and LM (solid line). (a) $N=30$; (b) $N=40$; (c) $N=50$. All cases use $1110$ training data points.}
\label{fig:ex02_optimizer_Adam-LBFGS-LM}
\end{center}
\end{figure}

\paragraph{\emph{Comparison of different optimizers.}}
The reasons why we choose Levenberg-Marquardt algorithm as our optimizer are two-fold. First, the LM algorithm is a combination of Gauss-Newton and gradient descent method which is suitable for nonlinear least squares problems. (The minimization of the loss function in the present paper is a nonlinear least square problem.) Meanwhile, the number of parameters to be trained in our proposed neural network is moderate (a few hundreds), so the cost per epoch for LM algorithm is acceptable. Second, the LM algorithm usually converges faster than commonly used optimizers such as Adam~\cite{KB14} and~L-BFGS \cite{LN89}. Here, we compare the training performance for three different optimizers (Adam, L-BFGS, LM) by showing the corresponding training loss evolutions in Figure~\ref{fig:ex02_optimizer_Adam-LBFGS-LM}. We use the previous setup and fix the number of training points $M=1110$ but vary the number of neurons from $N=30$ to $50$.  One can see that, the LM optimizer can effectively reduce the loss to $O(10^{-10})$ within $3000$ epochs when the number of neurons increases. In contrast, the Adam and L-BFGS optimizers reduce the loss values more slowly, and  barely achieve the losses of the magnitude $O(10^{-2})$ and $O(10^{-4})$ even up to $10^5$ epochs. Although not shown here, the final relative errors of LM algorithm show about three orders of magnitude smaller than the ones obtained by the Adam or L-BFGS.

\paragraph{\emph{Comparison of different augmented inputs.}}
In the third experiment, we demonstrate the robustness of present cusp-enforced level set function augmented input $\phi_a=\abs{\phi}$. Here, we keep $\eta=10$ but choose $\alpha=0$ and impose Dirichlet boundary condition on $\pd\Omega$ for simplicity. We also set  $\gamma=1$ so the flux jump $\dbblk{\beta\partial_n u}$ is zero while the solution $u$ still has discontinuous first derivatives to focus on the expressibility of the present network. We compare the relative errors and the losses of using either $\phi$ or $\phi_a$ as the augmented input in a fixed shallow neural network with the number of neurons $N=40$. The total training points used is $M=1110$ (or $M_0 =30$). The results are shown in Table~\ref{tab:ex02_ctsphi&absphi_eta10} where the used augmented input is listed in the first column. One can see that the prediction accuracy for the level set function input $\phi$ is quite poor. The relative errors for $\phi$ and $\phi_a$ input are $O(10^{-1})$ and $O(10^{-5})$, respectively, so the latter significantly outperforms the former. Therefore, the present cusp-enforced augmented feature input is indeed more accurate and capable of tackling the interface problem with discontinuous first derivatives.

\begin{table}[t]
\begin{center}
    \begin{tabular}{c|c|c|ccc}
    \hline
Augmented input	&  $\|u_\mathcal{N}-u\|_\infty/\|u\|_\infty$ & $\|u_\mathcal{N}-u\|_2/\|u\|_2$ & $\Loss(\bft)$ \\
    \hline
    $\phi$ 	   &  $8.01\times10^{-1}$ & $9.13\times10^{-1}$ & $5.57\times10^{-2}$ \\
    $\phi_a$   &  $2.98\times10^{-5}$ & $3.32\times10^{-5}$ & $9.17\times10^{-9}$ \\
   \hline
    \end{tabular}
\caption{Relative errors and training loss for the shallow network with different augmented inputs, $\phi$ and $\phi_a$. Here, $\alpha=0$, $\eta=10$, and $\gamma=1$ in Example~2. $(M_0, M)=(30, 1110)$,  $(L,N, N_\bft)=(1,40,200)$ }
\label{tab:ex02_ctsphi&absphi_eta10}
    \end{center}	
\end{table}

\begin{figure}[t]
\begin{center}
\includegraphics[scale=0.30]{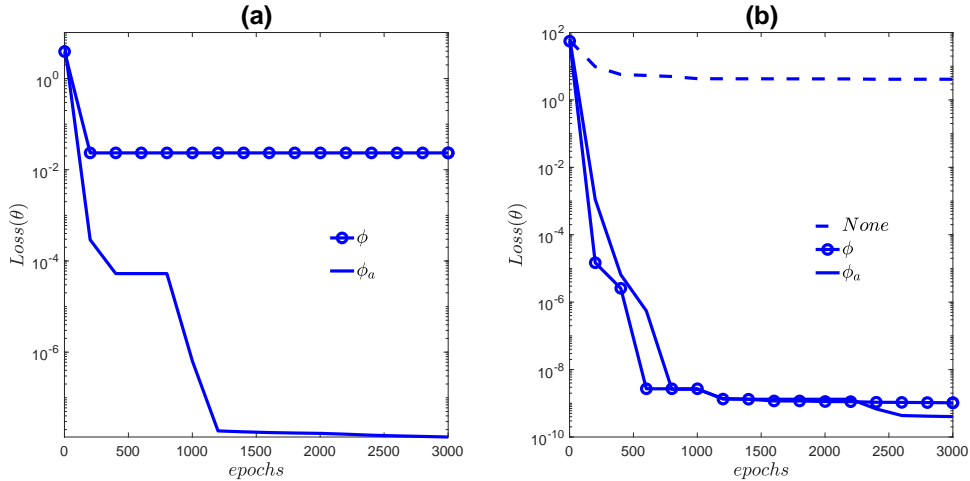}
\caption{(a)The evolutions of training $\Loss(\bft)$ corresponding to Table~\ref{tab:ex02_ctsphi&absphi_eta10}. (b) The evolutions of training $\Loss(\bft)$ corresponding to Table~\ref{tab:ex02_NN&ctsphi&absphi_eta1}. }
\label{fig:ex02_eta10&1_NN&Ctsphi&Absphi_loss}
\end{center}
\end{figure}

We also show the evolutionary plots of training loss for the two cases in Figure~\ref{fig:ex02_eta10&1_NN&Ctsphi&Absphi_loss}(a). After a few hundreds of epochs, the training loss for the case with augmented input $\phi$  becomes sluggish while the one with $\phi_a$ input continues to go down afterwards and reaches to the order of $10^{-8}$ eventually.

To further investigate  the power of function expressibility on the proposed cusp-enforced level set function augmentation, we consider a special case with $\eta = 1$ ($\beta^-=\beta^+=1$) and $\gamma = 1$ so that the jumps $\dbblk{\beta}=0$ and $\dbblk{\beta\partial_n u}=0$ simultaneously. One can immediately see from Eq.~(\ref{eq:jump_cond1}) that the normal derivative jump of $u$ equals to zero too, i.e., $\dbblk{\partial_n u}=0$. In this case, the solution $u$ is continuously differentiable across the interface $\Gamma$ so one might wonder if the level set function augmentation makes any differences. Table~\ref{tab:ex02_NN&ctsphi&absphi_eta1} shows the results for a shallow network with $\phi$, $\phi_a$ and  without augmented input (denoted by ``None"). For the one without augmented variable, the input is solely the position $\bx$. To have the same number of parameters used in the network, the one without augmented input uses $N=50$ neurons while the ones with augmented level set function input use $N=40$ neurons. Despite the fact that the solution is $C^1$, the network with solely $\bx$ input cannot train the solution properly as the training loss remains $O(1)$ (see Figure~\ref{fig:ex02_eta10&1_NN&Ctsphi&Absphi_loss}(b)) so the relative errors are greater than $5\%$. Again, the errors with $\phi_a$ augmented input are smaller than the ones with $\phi$ input in two orders of magnitude; that is, $O(10^{-5})$ versus $O(10^{-3})$. One can perceive that the network with cusp-enforced level set function augmentation still can predict the solution more accurately even though it is designed to capture the first-order derivatives correctly while the second-order derivatives are discontinuous across the interface in this example.

We also show the overall training time in the last column of Table~\ref{tab:ex02_NN&ctsphi&absphi_eta1}. Under the same setting, the training time per epoch using cusp-capturing PINN is indeed more costly than the one using the PINN (without any augmented input). However, as discussed earlier, if we use merely PINN, we are unable to train the network successfully even though the solution has the zero flux jump.

\begin{table}[h!]
\begin{center}
\begin{tabular}{c|c|c|c|c}
\hline
  Augmented input &   $\|u_\mathcal{N}-u\|_\infty / \|u\|_\infty$ & $\|u_\mathcal{N}-u\|_2 / \|u\|_2$ & $\Loss(\bft)$ & Elapsed time\\
    \hline
 None (PINN) \,   &  $1.59\times10^{-1}$ & $9.15\times10^{-2}$ & $4.34\times10^{0}$    & 10.9(s)\\
$\phi$     		  &  $8.55\times10^{-3}$ & $5.43\times10^{-3}$ & $9.07\times10^{-9}$   & 17.8(s)\\
$ \phi_a $ 		  &  $1.65\times10^{-5}$ & $1.91\times10^{-5}$ & $4.61\times10^{-10}$  & 22.1(s)\\
  \hline
\end{tabular}
\caption{Relative errors,  training losses, and total training time for the shallow network with an augmented input ($\phi$ or $\phi_a$) or without augmented input. Here,  $\alpha=0$, $\eta=1$, and $\gamma=1$ in Example 2. $(M_0, M)=(30, 1110)$.}
\label{tab:ex02_NN&ctsphi&absphi_eta1}
    \end{center}	
\end{table}

\begin{figure}[h]
\begin{center}
\includegraphics[scale = 0.32]{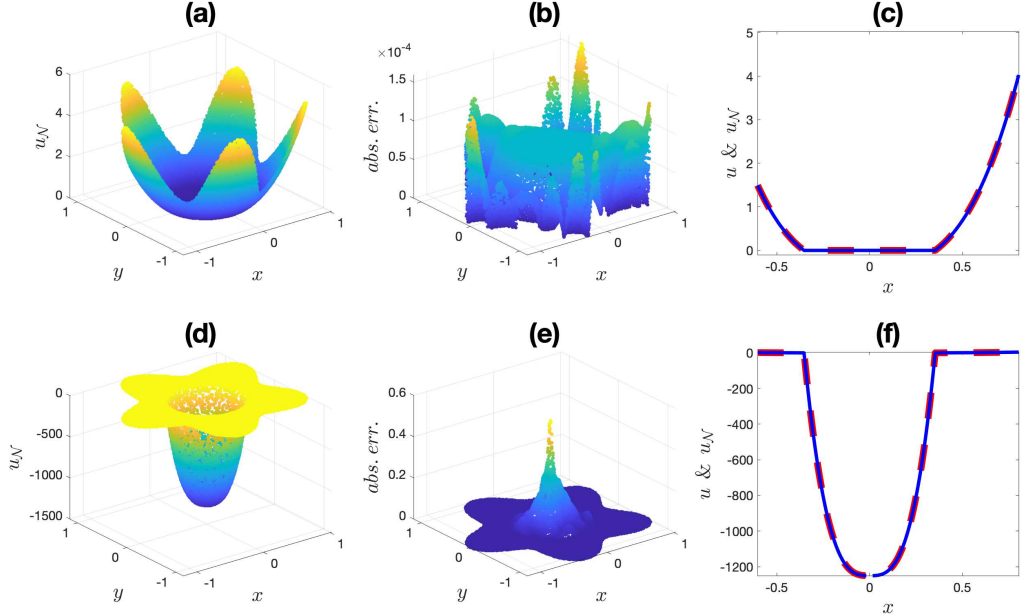}
\caption{(a) and (d): The profile of $u_\mathcal{N}$; (b) and (e): Absolute point-wise error $\abs{u_\mathcal{N}-u}$; (c) and (f): Cross-sectional view of $u_\mathcal{N}$ (blue solid line) and $u$ (red dashed line) along the line $y=x$. The upper panel is for $\eta = 10^{4}$ with $(L,N,N_\bft)=(2,15,315)$, and the lower panel is for $\eta = 10^{-4}$ with $(L,N,N_\bft)=(2,28,952)$ in Example~3.}
\label{fig:fig_ex3_tab04}
\end{center}
\end{figure}

\paragraph{Example 3.}
The third example illustrates that the present method is applicable for solving interface problems with high-contrast coefficients defined on irregular domains. We consider a five-fold flower region $\Omega=\{(x(r,\bft),\,y(r,\bft))\in\mathbb{R}^2\,|\,r(\bft)\le\,1-0.2\cos(5\bft)\}$ with an embedded interface, $\Gamma=\{(x,\,y)\in\mathbb{R}^2\:|\:x^2+y^2=\frac{1}{4}\}$. As in Example 2, the coefficient $\beta$ is defined in a piecewise-constant manner. The exact solution $u$ is defined as

\begin{equation}
u(x,y) = \left\{ \begin{array}{ll}   \frac{1}{\beta^-}\left(\left(x^2+y^2\right)^{\frac{3}{2}}-\frac{1}{8}\right), & (x,y)\in\Omega^-, \\ \frac{3}{\beta^+}\left(\left(x^2+y^2\right)^{\frac{3}{2}} - \frac{1}{8}\right), & (x,y)\in\Omega^+,\end{array}\right.
\end{equation}
and the Dirichlet boundary condition is imposed for simplicity. This problem was similarly studied by Wang et. al.~\cite{WZ20} using deep Ritz method on a square domain with an embedded circular interface. Again, the contrast ratio is defined by  $\eta=\beta^-/\beta^+$, and we fix $\beta^+=1$ so $\eta=\beta^-$. Here, we consider two high-contrast ratios; namely $\eta=10^{-4}$ and $10^{4}$. The cusp-enforced level set function is chosen as $\phi_a(x,y) = \abs{4(x^2+y^2)-1}$.  We generate $M = 1498$ training data ($M_I=1138$, $M_B=240$, and $M_\Gamma=120$) for the case of $\eta=10^4$, and employ the networks comprising from single to three hidden layers. The number of neurons for each network is chosen such that the number of trainable parameters $N_\bft$ is almost the same. As shown in Table~\ref{tab:ex3_highcontrastbeta}, for the contrast ratio $\eta = 10^{4}$, all network solutions can achieve accurate prediction with relative $L^2$ errors ranging from $O(10^{-4})$ to $O(10^{-5})$, which outperform the results obtained in \cite{WZ20}. However, for the contrast ratio $\eta = 10^{-4}$, the magnitude of exact solution $u$ in $\Omega^-$  is of the order $O(10^3)$ which is much larger than the solution in $\Omega^+$ of $O(1)$ (see also in Figure~\ref{fig:fig_ex3_tab04}(d)). So we have to use more neurons and training points ($M = 2959$ with $M_I=2519$, $M_B=240$, and $M_\Gamma=200$) to train the networks.  In this case, the relative errors range from  $O(10^{-3})$ to $O(10^{-4})$.  In addition, we depict the network solution profile, absolute point-wise error, and the cross-sectional view along $y = x$ in Fig.~\ref{fig:fig_ex3_tab04}. The upper and lower panels are for the contrast ratio $\eta = 10^{4}$ and $10^{-4}$, respectively. One can see that, without paying extra numerical efforts, the present model is able to tackle the interface problems in irregular domains thanks to the mesh-free advantage of neural network approximation. On the other hand, it could be quite tedious in implementation for traditional grid-based methods to handle such problems.

\begin{table}[h!]
\begin{center}
   \begin{tabular}{c|c|ccl}
   \hline
$\eta=\beta^-/\beta^+$ & $(L,N,N_\bft)$ & $\|u_\mathcal{N}-u\|_\infty/ \|u\|_\infty$ & $\|u_\mathcal{N}-u\|_2/ \|u\|_2$ & $\Loss(\bft)$ \\
    \hline
  			& $(1,63,315)$ & $3.29\times10^{-5}$ & $3.29\times10^{-5}$ & $1.35\times10^{-9}$ \\
$10^{4}$ 	& $(2,15,315)$ & $3.65\times10^{-5}$ & $3.82\times10^{-5}$ & $7.29\times10^{-10}$ \\
  			& $(3,11,319)$ & $8.73\times10^{-5}$ & $1.24\times10^{-4}$ & $3.32\times10^{-9}$ \\
  	\hline
 		  	& $(1,190,950)$ & $4.25\times10^{-3}$ & $1.42\times10^{-3}$ & $6.25\times10^{-9}$ \\
$10^{-4}$ 	& $(2,28,952)$  & $3.88\times10^{-4}$ & $1.28\times10^{-4}$ & $1.82\times10^{-11}$ \\
 		  	& $(3,20,940)$  & $1.50\times10^{-3}$ & $5.07\times10^{-4}$ & $3.05\times10^{-10}$ \\
    \hline
    \end{tabular}
    \caption{Relative errors and training losses in Example~3.}\label{tab:ex3_highcontrastbeta}
    \end{center}	
\end{table}

\paragraph{Example 4.}
In the fourth examples, we deal with the three-dimensional discontinuous variable-coefficient case and compare the accuracy of the present network solution with one of the immersed interface method (IIM) in~\cite{DIL03}. The domain is set as the cube $\Omega=[-1,1]\times[-1,1]\times[-1,1]$ in which the embedded interface is given by $\Gamma=\{(x,\,y,\,z)\in\mathbb{R}^3\:|\:x^2+y^2+z^2=r_0^2\}$. The exact solution $u$ and the variable-coefficient $\beta$ are chosen the same as in~\cite{DIL03},
\begin{equation}
u(x,y,z) = \left\{ \begin{array}{ll} r^2, & r<r_0, \\
r^2_0 + \frac{1}{b}\left(\frac{r^4}{2}+r^2-\frac{r^4_0}{2}-r^2_0\right), & r\ge r_0, \end{array}\right.\nonumber
\end{equation}
and
\begin{equation*}
\beta(x,y,z)=\left\{ \begin{array}{ll} r^2+1, & \quad r<r_0, \\ b, & \quad r\ge r_0, \end{array}\right.	
\end{equation*}
where $r_0=1/2$, $r=\sqrt{x^2+y^2+z^2}$, and the source term $f(x,y,z)=10r^2+6$. The solution satisfies the homogeneous jump conditions $\dbblk{u}=0$ and $\dbblk{\beta\partial_n u}=0$. However, the variable coefficient $\beta$ controlled by the parameter $b$ implies the discontinuity of the normal derivative $\dbblk{\partial_n u}$ at the interface $\Gamma$. The cusp-enforced level set function is chosen as $\phi_a(x,y,z) = \abs{4(x^2+y^2+z^2)-1}$.

For the sampling of training data points, we generate $M_I$ data points in the region $\Omega^+\cup\Omega^-$, and $M_B$ on the domain boundary ($M_B/6$ uniformly distributed training points on each face), while $M_\Gamma$ data points on the surface $\Gamma$ are generated by DistMesh~\cite{PS04}. In each of the following tests, the number of overall training points used is $M=3360$ ($M_I=800$, $M_B=2400$, and $M_\Gamma=160$).

Table~\ref{tab:ex04_varyeta_ShallowDeep} shows the relative errors and losses of the present method for three cases of $b=1$, $10$ and $1000$. Surprisingly, no matter how large the parameter $b$ is,  the present method with single- or multiple-hidden-layer structure gives accurate network predictions with the relative $L^\infty$ and $L^2$ errors of the magnitude $O(10^{-6})$. Here, we also present the results produced by IIM~\cite{DIL03} using $104\times104\times104$ uniformly distributed grid points. It should be noted that, in 3D IIM, the total number of degree of freedom (unknowns) is $104^3$ while the number of trainable parameters is just about $240$ for the present method. One can clearly see that our results outperform the ones obtained by IIM in almost two orders of magnitude in the relative $L^\infty$ error.

\begin{table}[h!]
\begin{center}
   \begin{tabular}{c|c|ccc}
   \hline
    $b$ & $(L, N, N_\bft)$ & $\|u_\mathcal{N}-u\|_\infty/\|u\|_\infty$ & $\|u_\mathcal{N}-u\|_2/\|u\|_2$ & $\Loss(\bft)$ \\
    \hline
 			& $(1, 40, 240)$ & $1.90\times10^{-6}$ & $2.17\times10^{-6}$ & $7.13\times10^{-11}$ \\
 	$1$ 	& $(2, 12, 228)$ & $1.15\times10^{-6}$ & $1.77\times10^{-6}$ & $7.53\times10^{-11}$ \\
 			& $(3, 9, 234)$  & $1.49\times10^{-6}$ & $1.52\times10^{-6}$ & $5.86\times10^{-11}$ \\
 			& IIM 			 & $9.59\times 10^{-5}$ & & \\
   	\hline
 			& $(1, 40, 240)$ & $2.48\times10^{-6}$ & $1.92\times10^{-6}$ & $3.68\times10^{-11}$ \\
 	$10$ 	& $(2, 12, 228)$ & $3.82\times10^{-6}$ & $1.84\times10^{-6}$ & $5.59\times10^{-11}$ \\
 			& $(3, 9, 234)$  & $3.95\times10^{-6}$ & $3.11\times10^{-6}$ & $3.46\times10^{-11}$ \\
 			& IIM 			 & $1.01\times 10^{-4}$ & & \\
   	\hline
 			& $(1, 40, 240)$ & $5.04\times10^{-6}$ & $3.51\times10^{-7}$ & $7.83\times10^{-11}$ \\
 	$1000$  & $(2, 12, 228)$ & $5.89\times10^{-6}$ & $6.38\times10^{-7}$ & $1.43\times10^{-10}$ \\
 			& $(3, 9, 234)$ 	 & $4.40\times10^{-6}$ & $5.95\times10^{-7}$ & $2.21\times10^{-10}$ \\
 			& IIM 			 & $1.61\times 10^{-4}$ & & \\
   \hline
    \end{tabular}
\caption{Relative errors and training losses in Example~4. The results produced by IIM use $104\times104\times104$ grid points.} \label{tab:ex04_varyeta_ShallowDeep}
    \end{center}	
\end{table}

\paragraph{Example 5.} In this example, we consider a problem of dimension $d=6$ to show that the present method is able to solve high-dimensional problems. Same problem was also solved in~\cite{LCLHL22} using a shallow Ritz method. Here we consider the domain $\Omega$ as a $6$-sphere of radius $0.6$ enclosing a smaller $6$-sphere of radius $0.5$ as $\Omega^-$. The cusp-enforced level set function is chosen as $\phi_a(\bx) = \abs{\left(\|\bx\|_2/0.5\right)^2 - 1}$, where $\bx = (x_1, x_2,x_3,x_4,x_5,x_6)$. We fix $\alpha=0$, a constant coefficient $\beta(\bx)=1$, and the exact solution is defined as
\begin{align}
    u(\mathbf{x}) =
    \left\{
    \begin{array}{ll}
        \exp(0.5^2-\|\mathbf{x}\|_2^2) + \sum^5_{i=1}\sin(x_i)& \mathbf{x} \in\Omega^+,\\
        1+2\sin(0.5^2 - \|\mathbf{x}\|_2^2) + \sum^5_{i=1}\sin(x_i) & \mathbf{x} \in \Omega^-.\\
    \end{array}\right.
\end{align}
The right-hand side functions can be obtained using Eqs.~(\ref{eq:varelliptic})-(\ref{eq:bdc}).

We use a shallow network ($L=1$) structure with $M=2628$ points to train the network.
The results are shown in Table~\ref{tab:ex6D_Shallow}. Using $40$ neurons in the hidden layer (and correspondingly $360$ trainable parameters), the relative $L^{\infty}$ and $L^2$ errors are in the order of $O(10^{-6})$ and $O(10^{-7})$, respectively. This example shows that the present method is applicable to solve high-dimensional elliptic interface problems.

\begin{table}[h!]
\begin{center}
   \begin{tabular}{c|ccc}
   \hline
    $(N, N_\theta)$ & $\|u_\mathcal{N}-u\|_\infty/\|u\|_\infty$ & $\|u_\mathcal{N}-u\|_2/\|u\|_2$ & $\Loss(\theta)$ \\
    \hline
 	$(10, 90)$   & $2.16\times10^{-3}$ & $1.37\times10^{-3}$ & $1.08\times10^{-4}$ \\
 	$(20, 180)$  & $7.69\times10^{-4}$ & $2.45\times10^{-4}$ & $1.48\times10^{-6}$ \\
 	$(30, 270)$  & $9.54\times10^{-5}$ & $3.79\times10^{-5}$ & $1.51\times10^{-8}$ \\
 	$(40, 360)$  & $1.86\times10^{-6}$ & $6.90\times10^{-7}$ & $5.77\times10^{-11}$ \\
 	\hline
    \end{tabular}
\caption{Relative errors and losses with $M=2628$ training data points where $(M_I, M_B, M_\Gamma)=(500, 1064, 1064)$ in Example~5.} \label{tab:ex6D_Shallow}
    \end{center}	
\end{table}

\paragraph{Example 6.} Next, we take an example in \cite{ASL20} that has its solution being discontinuous and make an accuracy comparison with the recent mesh-free methods ~\cite{ASL20, Oruc21}. We consider a two-dimensional computational domain $\Omega=[-1,1]\times[-1,1]$ with an embedded circular interface $\Gamma$ represented by the zero level set function  $\phi(x,y) = x^2+y^2-\left(\frac{2}{3}\right)^2$. The exact solution is expressed as
\beq
u(x,y) = \left\{\bary{lc} \sin(4\pi x)\sin(4\pi y)+7, & (x,y)\in\Omega^-, \\ 5\exp(-x^2-y^2), & (x,y)\in\Omega^+, \eary\right.\nonumber
\eeq
and the coefficient $\beta$ is a piecewise constant with $\beta^-=2$ and $\beta^+=3$ in $\Omega^-$ and $\Omega^+$, respectively.  This example also aims to demonstrate the applicability of the proposed method presented in the Remark since the above analytic solution is discontinuous across the interface.  Following the procedures in the Remark, we use a shallow network with $100$ neurons and $1000$ random points $\bx_\Gamma$ on the interface to train the function $V(\bx)$ satisfying $V(\bx_\Gamma)=-\dbblk{u}(\bx_\Gamma)$.  Once $V$ is available (thus $v$ is obtained), we apply the present cusp-capturing PINN to solve Eqs~(\ref{eq:varelliptic2})-(\ref{eq:bdc2}) to obtain the solution $w$. Then we can recover the solution $u=v+w$. Table~\ref{tab:discts2d_Oruc} presents the $L^\infty$ errors of the proposed method and two other non-neural network mesh-free methods, including the local mesh-free method based on LMM2P in \cite{ASL20} and the Pascal polynomials-based multiple-scale approach in \cite{Oruc21}.

\begin{figure}[h]
\begin{center}
\includegraphics[width=\textwidth]{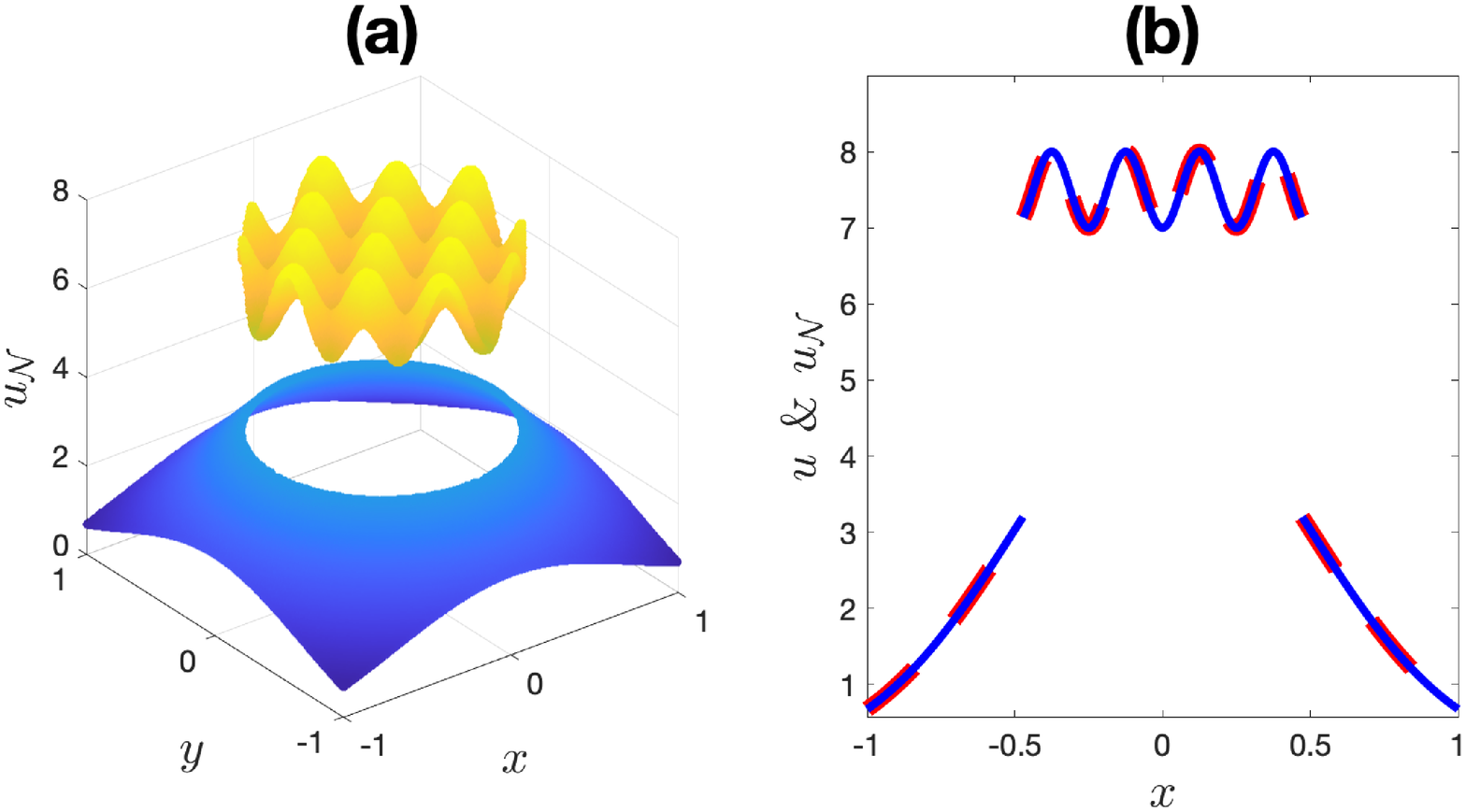}
\caption{The corresponding network solution plots in Example 6. $(L,N,N_\theta)=(3,20,940)$, $(M_I, M_B, M_\Gamma)=(3635, 600, 300)$. (a) The profile of the network solution $u_\mathcal{N}$; (b) Cross-sectional view of $u_\mathcal{N}$ (blue solid line) and $u$ (red dashed line) along the line $y=x$.}
\label{fig:fig_ex6_D20L3}
\end{center}
\end{figure}

The present networks with the number of hidden layer $L=1, 2$  use exactly same number of trainable parameters $N_\theta=775$ and same number of total training points $M=3150$  ($M_I=2550$, $M_B=400$, and $M_\Gamma=200$) which give the $L^\infty$ errors ranging from the magnitude $O(10^{-3})$ to $O(10^{-4})$. Here,  using a deeper network seems to predict more accurate results than the shallow one under the same number of trainable parameters used. Therefore, we use a deep network $(L,N,N_\theta)=(3,20,940)$ with $M=4535$ ($M_I=3635$, $M_B=600$, and $M_\Gamma=300$) training points to reduce the $L^\infty$ error to the magnitude of $O(10^{-5})$, where the solution profile $u_\mathcal{N}$ and its cross-sectional view along the line $y=x$ are shown in Figure~\ref{fig:fig_ex6_D20L3}. In the Table, the number of nodes indicates the number of mesh-free points used in these methods which works like the number of training points $M$ used in the present method.  One can immediately see that our numerical results are slightly more accurate than the ones in \cite{ASL20} and less accurate than the ones obtained in \cite{Oruc21}.

\begin{table}[h!]
\begin{center}
   \begin{tabular}{c|c|cc|cc}
   \hline
    $(L, N, N_\theta,M)$ & Present & No. nodes &  Ahmad et al.~\cite{ASL20} & No. nodes & Oru\c{c}~\cite{Oruc21} \\
    \hline
   $(1,155, 775, 3150)$  & $4.31\times10^{-3}$ & $1600$ &  $8.75\times10^{-3}$ & $1365\,(23)$ & $1.54\times10^{-3}$ \\
   $(2, 25, 775, 3150)$  & $3.78\times10^{-4}$ & $6400$ &  $1.52\times10^{-3}$ & $2490\,(25)$ & $1.04\times10^{-4}$ \\
   $(3, 20, 940, 4535)$  & $3.72\times10^{-5}$ & $25600$ &  -                   & $4065\,(27)$ & $4.08\times10^{-6}$ \\
 	\hline
    \end{tabular}
\caption{Comparison of $L^\infty$ errors using the present method and two recent mesh-free methods \cite{ASL20, Oruc21} in Example~6. The number in the parentheses represents the number of $m$ with the highest degree of polynomial $m-1$ used in~\cite{Oruc21}.} \label{tab:discts2d_Oruc}
    \end{center}	
\end{table}

\begin{figure}[!]
\begin{center}
\includegraphics[width=\textwidth]{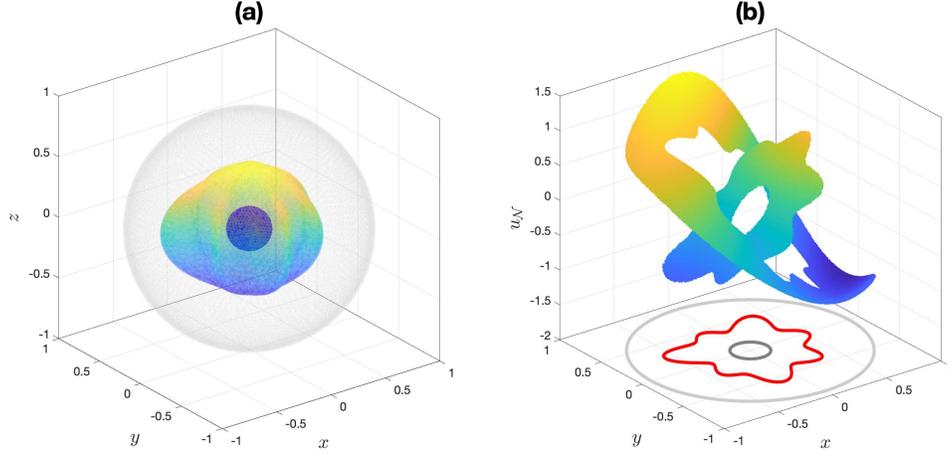}
\caption{(a) The illustration of domain and interface geometries in Example~7. (b) The cross-sectional view of the network solution $u_\mathcal{N}$ on $z=0$. The red and grey curves indicate the corresponding cross-sectional interface and domain boundaries, respectively.}
\label{fig:fig_ex7_S100M2460_twoplots}
\end{center}
\end{figure}

\paragraph{Example 7.} The last example is taken from~\cite{BG20}, in which we consider a spherical shell $\Omega=\{(x,y,z)\in\mathbb{R}^3\:|\: 0.151^2\le x^2+y^2+z^2\le 0.911^2 \}$ where a complex embedded interface $\Gamma$ is represented by the zero level set of the level set function
\beq
\phi(x,y,z) = \sqrt{x^2+y^2+z^2}-r_0\left(1+\left(\frac{x^2+y^2}{x^2+y^2+z^2}\right)^2\sum^3_{k=1}a_k\cos\left(n_k\left(\tan^{-1}\left(\frac{y}{x}\right)-\theta_k\right)\right)\right), \nonumber
\eeq
and the setup of parameters is shown as follows:
\beq
r_0=0.483,\:\: \left(\bary{c}a_1 \\ a_2 \\ a_3 \eary\right)=\left(\bary{c}0.1 \\ -0.1 \\ 0.15 \eary\right), \:\: \left(\bary{c} n_1 \\ n_2 \\ n_3 \eary\right)=\left(\bary{c}3 \\ 4 \\ 7 \eary\right),\:\mbox{and}\: \left(\bary{c}\theta_1 \\ \theta_2 \\ \theta_3 \eary\right)=\left(\bary{c} 0.5 \\ 1.8 \\ 0 \eary\right). \nonumber
\eeq
The illustration of domain and interface geometry can be found in Figure~\ref{fig:fig_ex7_S100M2460_twoplots}(a). Note that, the dark-shading region enclosed in the interface is the inner boundary of the domain $\Omega$.
 We choose the following same solution $u$ and the coefficient $\beta$ as in \cite{BG20}
\beqs
u(x, y, z) &=& \left\{ \bary{ll} \sin(2x)\cos(2y)e^z, & (x,y,z)\in\Omega^-, \\
\left(16\left(\frac{(y-x)}{3}\right)^5-20\left(\frac{(y-x)}{3}\right)^3+5\left(\frac{(y-x)}{3}\right)\right)\log(x+y+3)\cos(z), & (x,y,z)\in\Omega^+, \eary\right.\nonumber\\
\beta(x,y,z)&=&\left\{\bary{ll} 10\left(1+\frac{1}{5}\cos\left(2\pi(x+y)\right)\sin\left(2\pi(x-y)\right)\cos(z)\right), & (x,y,z)\in\Omega^-, \\ 1, & (x,y,z)\in\Omega^+. \eary\right.\nonumber
\eeqs
The right-hand side functions can be obtained using Eqs.~(\ref{eq:varelliptic})-(\ref{eq:bdc}). 

Again, the above analytic solution is obviously discontinuous across the interface so we have to follow the solution procedures discussed in the Remark to obtain the approximate network solution. First, we use a shallow network with $100$ neurons and $752$ training points (generated by DistMesh~\cite{PS04}) $\bx_\Gamma$ on the interface to train the function $V(\bx)$ satisfying $V(\bx_\Gamma)=-\dbblk{u}(\bx_\Gamma)$.  Once $V$ is available (thus $v$ is obtained), we solve Eqs~(\ref{eq:varelliptic2})-(\ref{eq:bdc2}) by applying the present cusp-capturing PINN with one-hidden-layer and $2460$ training data points ($M_I=801$, $M_B=907$, and $M_\Gamma=752$) to train the solution $w$. After that, we obtain the network approximate solution $u=v+w$. Table~\ref{tab:ex7_complex_interface3d} shows the relative $L^\infty$ and $L^2$ errors for different number of neurons $N$ used in the hidden layer. One can see that, using merely $25$ neurons in the hidden layer (correspondingly $150$ trainable parameters), the relative errors and training losses are of the magnitudes $O(10^{-4})$ and $O(10^{-7})$, respectively. The relative errors can be reduced to the magnitude $O(10^{-6})$ when the number of neurons increases to $100$. Figure~\ref{fig:fig_ex7_S100M2460_twoplots}(b) shows the cross-sectional profile of the network solution on the hyperplane $z=0$. As a result, the present method is indeed applicable for solving elliptic interface problems in irregular domain with complex interface subject to nonzero solution jump condition.

\begin{table}[h!]
\begin{center}
   \begin{tabular}{c|ccc}
   \hline
	$(N,N_\bft)$  & $\|u_\mathcal{N}-u\|_\infty/ \|u\|_\infty$ & $\|u_\mathcal{N}-u\|_2/ \|u\|_2$ & $\Loss(\bft)$ \\
    \hline
	$(25,150)$   & $7.55\times10^{-4}$ & $5.03\times10^{-4}$ & $1.16\times10^{-7}$ \\
	$(50,300)$   & $3.04\times10^{-5}$ & $1.09\times10^{-5}$ & $1.54\times10^{-9}$ \\
	$(100,600)$  & $3.07\times10^{-6}$ & $1.05\times10^{-6}$ & $1.36\times10^{-11}$ \\
    \hline
    \end{tabular}
    \caption{Relative errors and training losses in Example~7.} \label{tab:ex7_complex_interface3d}
    \end{center}	
\end{table}

\section{Conclusion}
\label{sec:conclusion}
We propose a cusp-capturing physics-informed neural network for solving the discontinuous-coefficient elliptic interface problems. By introducing a cusp-enforced level set function as an additional feature input to the network, the predicted solution by the network can retain the inherent properties of the solution which is continuous but the normal derivative has a jump discontinuity on the interface. The training procedure uses the LM-based optimizer to minimize the loss function comprising mean squared errors of the equation residual, the interface condition, and the boundary condition in the same spirit as the physics-informed neural networks. We conduct a series of numerical tests to show the accuracy of the present network, with particular emphasis on the number of neurons and training points, and the effectiveness of the cusp-capturing technique. A high-contrast coefficient interface problem is included in our numerical experiments, and the accuracy outperforms the one obtained in previous work. The present network is efficient in terms of network structure since one hidden layer with a moderate number of neurons and sufficiently enough training data points can achieve quite accurate predictions.   The results are also comparable to traditional grid-based methods, such as the immersed interface method. Besides, if the solution is discontinuous across the interface, we can simply incorporate an additional supervised learning task for solution jump approximation into the present network without much difficulty. In the future, we shall apply the present network method to practical applications where traditional grid-based methods are difficult to implement and extend to the time-dependent discontinuous-coefficient interface problems. Meanwhile, using functions other than level sets to represent interfaces for handling the $C^0$-interfaces and considering multiple interfaces is beyond the scope of this paper and is certainly worthy exploring in the future.

\section*{Acknowledgement}
Y.-H. Tseng, T.-S. Lin, W.-F. Hu, and M.-C. Lai acknowledge the supports by National Science and Technology Council, Taiwan, under research grants 111-2115-M-390-002, 111-2628-M-A49-008-MY4, 111-2115-M-008-009-MY3, and 110-2115-M-A49-011-MY3, respectively. T.-S. Lin and W.-F. Hu also acknowledge the supports by National Center for Theoretical Sciences, Taiwan.

\end{document}